\documentclass[12pt]{amsart}
\usepackage{amssymb,amsthm,amsmath,epsfig,latexsym}
\usepackage{calc}

\usepackage{ifthen}
\newcommand{\showcomments}{yes}
\newsavebox{\commentbox}
\newenvironment{comment}%
{\ifthenelse{\equal{\showcomments}{yes}}%
{\footnotemark
    \begin{lrbox}{\commentbox}
    \begin{minipage}[t]{1.25in}\raggedright\sffamily\upshape\tiny
    \footnotemark[\arabic{footnote}]}
{\begin{lrbox}{\commentbox}}}%
{\ifthenelse{\equal{\showcomments}{yes}}%
{\end{minipage}\end{lrbox}\marginpar{\usebox{\commentbox}}}
{\end{lrbox}}}

\begin{document}

\newcommand{\mmbox}[1]{\mbox{${#1}$}}
\newcommand{\proj}[1]{\mmbox{{\mathbb P}^{#1}}}
\newcommand{\affine}[1]{\mmbox{{\mathbb A}^{#1}}}
\newcommand{\Ann}[1]{\mmbox{{\rm Ann}({#1})}}
\newcommand{\caps}[3]{\mmbox{{#1}_{#2} \cap \ldots \cap {#1}_{#3}}}
\newcommand{\N}{{\mathbb N}}
\newcommand{\Z}{{\mathbb Z}}
\newcommand{\R}{{\mathbb R}}
\newcommand{\Tor}{\mathop{\rm Tor}\nolimits}
\newcommand{\Ext}{\mathop{\rm Ext}\nolimits}
\newcommand{\Hom}{\mathop{\rm Hom}\nolimits}
\newcommand{\im}{\mathop{\rm Im}\nolimits}
\newcommand{\rank}{\mathop{\rm rank}\nolimits}
\newcommand{\supp}{\mathop{\rm supp}\nolimits}
\newcommand{\arrow}[1]{\stackrel{#1}{\longrightarrow}}
\newcommand{\CB}{Cayley-Bacharach}
\newcommand{\coker}{\mathop{\rm coker}\nolimits}
\sloppy
\newtheorem{defn0}{Definition}[section]
\newtheorem{prop0}[defn0]{Proposition}
\newtheorem{conj0}[defn0]{Conjecture}
\newtheorem{thm0}[defn0]{Theorem}
\newtheorem{lem0}[defn0]{Lemma}
\newtheorem{corollary0}[defn0]{Corollary}
\newtheorem{example0}[defn0]{Example}

\newenvironment{defn}{\begin{defn0}}{\end{defn0}}
\newenvironment{prop}{\begin{prop0}}{\end{prop0}}
\newenvironment{conj}{\begin{conj0}}{\end{conj0}}
\newenvironment{thm}{\begin{thm0}}{\end{thm0}}
\newenvironment{lem}{\begin{lem0}}{\end{lem0}}
\newenvironment{cor}{\begin{corollary0}}{\end{corollary0}}
\newenvironment{exm}{\begin{example0}\rm}{\end{example0}}

\newcommand{\defref}[1]{Definition~\ref{#1}}
\newcommand{\propref}[1]{Proposition~\ref{#1}}
\newcommand{\thmref}[1]{Theorem~\ref{#1}}
\newcommand{\lemref}[1]{Lemma~\ref{#1}}
\newcommand{\corref}[1]{Corollary~\ref{#1}}
\newcommand{\exref}[1]{Example~\ref{#1}}
\newcommand{\secref}[1]{Section~\ref{#1}}

\newcommand{\std}{Gr\"{o}bner}
\newcommand{\jq}{J_{Q}}



\title[Cayley-Bacharach and Evaluation Codes]{Cayley-Bacharach and Evaluation Codes on Complete Intersections}

\author{Leah Gold}
\address{Gold: Mathematics Department \\ Texas A\&M University \\
  College Station \\ TX 77843-3368 \\ USA}
\thanks{Gold is partially supported by an NSF-VIGRE postdoctoral fellowship}
\email{lgold@math.tamu.edu}
 
\author{John Little}
\address{Little: Department of Mathematics and Computer Science \\
College of the Holy Cross \\ Worcester \\ MA 01610 \\ USA}
\email{little@mathcs.holycross.edu}
 
\author{Hal Schenck}
\thanks{Schenck is partially supported by NSF Grant DMS 03--11142}
\address{Schenck: Mathematics Department \\ Texas A\&M University \\
  College Station \\ TX 77843-3368\\ USA}
\email{schenck@math.tamu.edu}
 

\begin{abstract}
\noindent
In \cite{jh1}, J. Hansen uses cohomological methods to find a lower
bound for the minimum distance of an evaluation code determined by a
reduced complete intersection in $\mathbb{P}^2$. In this paper,
we generalize Hansen's results from $\mathbb{P}^2$ to
$\mathbb{P}^m$; we also show that the hypotheses of \cite{jh1} may
be weakened.  The proof is succinct and follows by combining 
the Cayley-Bacharach Theorem and the bounds on evaluation codes obtained
in \cite{jh}.
\end{abstract}
\maketitle


\section{Introduction}\label{sec:one}
In \cite{drt}, Duursma, Renter\'\i a, and Tapia-Recillas compute the
block length and dimension of the Reed-Muller (or evaluation) code
determined by a zero-dimensional complete intersection $\Gamma
\subset \proj{m}$.  The words of the code $C(\Gamma)_a$ are obtained
by evaluating homogeneous polynomials of
degree $a$ at the points of $\Gamma$. When $\Gamma$ is determined by
two polynomials of degrees $d_1$, $d_2$ in $R=\mathbb{K}[x,y,z]$, Hansen
\cite{jh1} obtains a lower bound for the minimum distance of the code.
In particular, if $d_i \ge 3$ and
$$\max\{d_1-2, d_2-2\} \; \le \; a \; \le \; d_1+d_2-3,$$
then the code
$C(\Gamma)_a$ has minimum distance $d \ge d_1+d_2-a-1$.  The key point
is the observation that when one evaluates polynomials of degree $a$
with $a \le d_1+d_2-3,$ then the resulting evaluation vectors will be
linearly dependent.  In algebraic-geometric terms, this reflects the
fact that the points of $\Gamma$ fail to impose independent conditions on
polynomials of degree $a$. It turns out that this failure gives one
some room to correct transmission errors.

The main theme of this paper is that using the modern Cayley-Bacharach
Theorem due to Davis, Geramita, and Orecchia \cite{dgo} streamlines
the proof in \cite{jh1} substantially, and makes it easy to generalize
the results from $\proj{2}$ to $\proj{m}$.  In the $m=2$ case
the Cayley-Bacharach Theorem also allows us to drop the hypotheses
$\max\{d_1-2, d_2-2\} \le a$ and $d_i \ge 3$ of \cite{jh1}, so, in
particular, our result applies to Reed-Solomon codes.  We start off
with a quick review of evaluation codes, and a discussion of residual
schemes and the Cayley-Bacharach Theorem.

\subsection{Background on evaluation codes}
Let $V$ be a variety in $\proj{m}$ defined over the finite field
$\mathbb{F}_q$, with $\Gamma=\{p_1,\ldots,p_n\}$ a set of 
$\mathbb{F}_q$-rational points on $V$. Let
$R=\mathbb{F}_q[x_0,\ldots,x_m]$, and let $R_a$ denote the 
vector space of homogeneous polynomials of degree $a$.
Choose a degree $a$ and $f_0 \in R_a$
such that $f_0(p_i) \ne 0$ for all $i \in \{1,\ldots,n \}$.  The
evaluation map $e_a(\Gamma)$ is defined to be the linear map
\begin{eqnarray*}
e_a(\Gamma) : R_a & \to &\mathbb{F}_q^n\\
      f &\mapsto& \left(\frac{f(p_1)}{f_0(p_1)},\ldots,\frac{f(p_n)}{f_0(p_n)}\right)\\
\end{eqnarray*}
The image of $e_a(\Gamma)$ is a linear code of block length $n$, which
we will denote $C(\Gamma)_a$.  The codes $C(\Gamma)_a$ are called {\it
  evaluation codes} associated to $\Gamma$.  The minimum distance of
$C(\Gamma)_a$ is
$$d = d(C(\Gamma)_a) = \min_{w_1\ne w_2 \in C(\Gamma)_a} |w_1 - w_2|,$$
where $| \cdot |$ denotes the norm corresponding
to the {\it Hamming distance}, that is, the
number of nonzero entries in a word. 
Since $C(\Gamma)_a$ is closed under sums, the minimum
distance is also equal to the minimum over all nonzero codewords 
of the number of nonzero entries, or equivalently, the length of 
the words minus the
largest number of zero entries in any nonzero codeword.  

The Singleton bound implies that the minimum distance $d$, the 
block length $n$, and the dimension $k$ of a linear code satisfy
$d \le n - k + 1$.  Codes for which the upper bound are 
achieved are known as maximum distance separable, or MDS, codes.

\subsection{Background on the Cayley-Bacharach Theorem}
Let $\mathbb{K}$ be a field and 
suppose $\Gamma=\{p_1,\ldots, p_n\}$ is a set of distinct points in
$\mathbb{P}^m_\mathbb{K}$. As above, let $e_a$ be the evaluation map 
from the vector space $R_a$ of homogeneous polynomials of
degree $a$ to $\mathbb{K}^n$.  The kernel of this map consists of polynomials of 
degree $a$ which vanish on $\Gamma$, so the kernel is simply the degree 
$a$ piece of the ideal $I_\Gamma$. Hence we have an exact sequence of 
vector spaces
$$0 \longrightarrow (I_\Gamma)_a \longrightarrow R_a \stackrel{e_a}
{\longrightarrow} \mathbb{K}^n \longrightarrow {\rm coker}(e_a) \longrightarrow 0.$$
Using sheaf cohomology and
writing $\mathcal{I}_\Gamma$ for the sheaf of ideals corresponding to
$I_\Gamma$, we can identify ${\rm coker}(e_a) \cong H^1(\mathcal{I}_\Gamma(a))$.
Similarly, the kernel of $e_a$ can be identified with
$H^0(\mathcal{I}_\Gamma(a))$.  We will write $h^0(\mathcal{I}_\Gamma(a))$ to
denote the dimension of the kernel of $e_a$ as a vector space over
$\mathbb{K}$.  In similar fashion, the dimension of the vector space
$H^1(\mathcal{I}_\Gamma(a))$ will be denoted by $h^1(\mathcal{I}_\Gamma(a))$.
The set of points $\Gamma$ is said to {\em impose independent
  conditions} on polynomials of degree $a$ if
the rank of $e_a$ is $n$, that is, if $\dim \coker(e_a) =
h^1(\mathcal{I}_\Gamma(a))=0$.

The classical Cayley-Bacharach Theorem deals with the following
situation. Suppose that $Y_1$, $Y_2 \subset \mathbb{P}^2$ are plane
curves of degree $d_1$ and $d_2$ which intersect in a set $\Gamma$ of
$d_1d_2$ distinct points. Write $\Gamma = \Gamma' \cup \Gamma''$ with
$\Gamma'$ and $\Gamma''$ disjoint. If $a \le d_1+d_2-3$ is a
nonnegative integer, then the classical Cayley-Bacharach Theorem
asserts that the dimension of the vector space
$(I_{\Gamma'})_a/(I_{\Gamma})_a$ is equal to
$h^1(\mathcal{I}_{\Gamma''}(d_1+d_2-3-a))$, a measure of the failure
of $\Gamma''$ to impose independent conditions in degree
$d_1+d_2-3-a$.  For instance, if $d_1 = d_2 = 3$, $a=3$, and $\Gamma =
\Gamma' \cup \Gamma''$, with $\deg(\Gamma') = 8$ and $\deg(\Gamma'') =
1$, then the classical Cayley-Bacharach Theorem says that $\dim
(I_{\Gamma'})_a/(I_{\Gamma})_a =
h^1(\mathcal{I}_{\Gamma''}(0))$. Since $
h^1(\mathcal{I}_{\Gamma''}(0)) = 0$, every cubic that vanishes at the
8 points in $\Gamma'$ also vanishes at the point in $\Gamma''$.

To formulate the modern version of the Cayley-Bacharach Theorem, we
need to use the language of schemes. For background on schemes we
refer the reader to \cite{eh}, and for a thorough discussion of the
Cayley-Bacharach Theorem we recommend \cite{egh1}.

\begin{defn}[Residual schemes \cite{egh1}]\label{residual}
  Let $\Gamma$ be a zero-dimensional scheme with coordinate ring
  $A(\Gamma)$. Let $\Gamma' \subset \Gamma$ be a closed subscheme and
  $I_{\Gamma'}\subset A(\Gamma)$ be its ideal. The subscheme of
  $\Gamma$ {\em residual} to $\Gamma'$ is the subscheme defined by the
  ideal $$I_{\Gamma''}=Ann(I_{\Gamma'}/I_{\Gamma}).$$
\end{defn}

When $\Gamma$ is a complete intersection, $\Gamma'$ is residual to
$\Gamma''$ in $\Gamma$ iff $\Gamma''$ is residual to $\Gamma'$ in
$\Gamma$ (this need not be the case in general).  We are now ready to
state the version of the Cayley-Bacharach Theorem that we will use to
extend the minimum distance bound.


\begin{thm}[Davis-Geramita-Orecchia, \cite{dgo}]\label{modernCB}
  Let $\Gamma \subset \proj{m}$ be a complete intersection of
  hypersurfaces $X_1, X_2, \ldots, X_m$ of degrees $d_1, d_2, \ldots,
  d_m$ respectively, and let $\Gamma'$, $\Gamma'' \subset \Gamma$ be
  closed subschemes residual to one another. Set
  \[
  s = \left( \sum_{i=1}^m d_i \right) - m -1.
  \]
  Then, for any $a \geq 0$, we have
  \[
  h^0({\mathcal I}_{\Gamma'}(a)) - h^0({\mathcal I}_{\Gamma}(a)) = h^1({\mathcal
    I}_{\Gamma''}(s-a)).
  \]
\end{thm}
In \cite{dgo}, this theorem is proved with the assumption that the
ground field is infinite. When  $\Gamma$ is composed of 
$\mathbb{F}_q$-rational points, the statement holds by
interpreting the dimensions over $\overline{\mathbb{F}_q}$. 
If we use the monomial basis for $R_a$, then it is easy to 
see that the matrix of the evaluation map 
$e_a$ : $R_a \rightarrow \overline{\mathbb{F}_q}^n $ 
has entries in $\mathbb{F}_q$, so the dimensions of  the kernel 
and cokernel will be the same whether we work over the infinite field 
$\overline{\mathbb{F}_q}$ or the finite field  $\mathbb{F}_q$.
\section{Review of $\mathbb{P}^2$ result}\label{sec:two}
Let $\Gamma \subset \mathbb{P}^2$ be a reduced complete intersection
of two curves of degrees $d_1, d_2$ defined over $\mathbb{F}_q$.
Theorem 4.4 of \cite{jh1} tells us that if $d_i \ge 3$ and $\max \{d_i
-2 \} \le a \le d_1 + d_2 - 3,$ then the evaluation code $C(\Gamma)_a$
has minimum distance $d \ge d_1 + d_2 - a - 1$. The proof in
\cite{jh1} uses Serre duality to compute the dimension of a certain
cohomology group, which is why the hypothesis $a \ge \max
\{d_i -2 \}$ is needed; also useful is the following lemma 
(2.6 of \cite{jh1}):
\begin{lem}\label{Hlemma1}
  Let $\Gamma$ be a finite set of points in $\mathbb{P}^m$, with $|\Gamma|=\deg
  \Gamma$. Then for $j \geq |\Gamma|-1$, $h^1({\mathcal I}_{\Gamma}(j))=0.$
\end{lem}
What Hansen actually shows in the proof of Theorem 4.4 in \cite{jh1} is
that if $\Gamma \subseteq \mathbb{P}^2$ is a $(d_1$,$d_2)$ complete 
intersection, and $\Gamma' \subset \Gamma$ satisfies
$$|\Gamma'| \ge d_1d_2 - d_1 - d_2 + a + 4,$$
then the projection map
$\pi : C(\Gamma)_a \to C(\Gamma')_a$, obtained by deleting the
components of the codewords of $C(\Gamma)_a$ corresponding to the
points in $\Gamma''$, is injective. We warm up by using the
Cayley-Bacharach Theorem to give a slight improvement.
\begin{lem}\label{Hlemma2}
If $\Gamma' \subset \Gamma$ satisfies
$$|\Gamma'| \ge d_1d_2 - d_1 - d_2 + a + 2,$$
then the projection map
$\pi : C(\Gamma)_a \to C(\Gamma')_a$, obtained by deleting the
components of the codewords of $C(\Gamma)_a$ corresponding to the
points in $\Gamma''$, is injective.
\end{lem}
\begin{proof}
  Since $\Gamma$ is reduced, $|\Gamma| = d_1d_2$. Let $s = d_1+d_2-3$
  and let $\Gamma'$ be any subset of the points of $\Gamma$ such that
  $|\Gamma'| \geq d_1d_2 -s+a-1$. Then letting $\Gamma'' = \Gamma
  \setminus \Gamma'$ be the subscheme residual to $\Gamma'$ we have
  $$|\Gamma''| \leq d_1d_2 -(d_1d_2-s+a-1) = s-a+1.$$
  Since $s-a \geq
  |\Gamma''|-1$, Lemma~\ref{Hlemma1} tells us that $\Gamma''$ imposes
  independent conditions in degree $s-a$, so $h^1({\mathcal
    I}_{\Gamma''}(s-a)) = 0$. On the other hand, $\Gamma'$ and
  $\Gamma''$ are closed subschemes of $\Gamma$ residual to one
  another, so by Theorem~\ref{modernCB} we know that for any $a \geq
  0$,
  $$
  h^0({\mathcal I}_{\Gamma'}(a)) - h^0({\mathcal I}_{\Gamma}(a)) =
  h^1({\mathcal I}_{\Gamma''}(s-a)).$$
  The right hand side is zero, so
  $h^0({\mathcal I}_{\Gamma'}(a)) = h^0({\mathcal I}_{\Gamma}(a))$.
  In other words, $H^0({\mathcal I}_{\Gamma'}(a)) \simeq H^0({\mathcal
    I}_{\Gamma}(a))$, that is, $(I_{\Gamma'})_a = (I_{\Gamma})_a$.
  Hence the projection map $C(\Gamma)_a
  \stackrel{\pi}{\longrightarrow} C(\Gamma')_a$ is injective.
  Moreover, the map is injective for {\it all ways} of splitting
  $\Gamma$ as a union of $\Gamma'$ and $\Gamma''$ with the same
  cardinality as above.
\end{proof}

We claim that the result
$d \ge s - a + 2$ on the minimum distance now follows. To see
this, consider the case $|\Gamma'| = d_1 d_2-s+a-1$
and $|\Gamma''| = s-a+1$. Let $0 \ne f \in R_a$.  If $f$ is
nonzero at $s-a+2$ or more points in $\Gamma'$,
then we are done, so we assume that $f$ is only nonzero at $t$ points
in $\Gamma'$ with
$$1 \le t \le s - a + 1.$$
It suffices to see that $f$ must be
nonzero at $\ge s-a+2-t$ points in $\Gamma''$.  If
not, then $f$ is nonzero at $\le s-a+1-t$ points in
$\Gamma''$, so $f$ vanishes at $\ge |\Gamma''| -
(s-a+1-t) = t$ points of $\Gamma''$.  Then we can
subdivide $\Gamma$ into two new 0-cycles $\overline{\Gamma}'$ and
$\overline{\Gamma}''$ by exchanging $t$ points from $\Gamma''$
where $f$ vanishes with $t$ points from $\Gamma'$ where $f$ is
nonzero.  We obtain a new
decomposition $\Gamma = \overline{\Gamma}' \cup \overline{\Gamma}''$
such that $f$ vanishes at all the points in $\overline{\Gamma}'$.
From the previous proof, we know that $C(\Gamma)_a
\stackrel{\pi}{\longrightarrow} C(\overline{\Gamma}')_a$ is injective,
so $f$ must vanish on all of $\Gamma$.  It follows that
$d \ge s - a + 2$. If $|\Gamma'| > d_1 d_2-s+a-1$, then we
can apply the same argument to any subset of $\Gamma'$ of 
size $d_1 d_2-s+a-1$ to obtain the bound.

\section{Main theorem}\label{sec:three}
We are now ready to prove the main result of this paper: Hansen's
bound generalizes to reduced complete intersections in $\mathbb{P}^m$.
This can be proved along the lines just sketched for the
$\mathbb{P}^2$ case.  However, the proof is shorter if we utilize the
criteria of \cite{jh} (Proposition 6 and Theorem 8). In the 
language of this paper, the result is:

\begin{prop}\label{Kprop} 
Let $\Gamma$ be a subset of points in $\mathbb{P}^m$, and let
  $C(\Gamma)_a$ be the evaluation code defined in \S 1.
  For $i \ge 1$, $d(C(\Gamma)_a) \ge deg(\Gamma)-i+1$
  iff $h^0(\mathcal{I}_\Gamma(a))=h^0(\mathcal{I}_{\Gamma'}(a))$ for all
  $\Gamma' \subset \Gamma$ with $|\Gamma'|=i$. Furthermore, $C(\Gamma)_a$ 
  is an MDS
  code iff $h^0(\mathcal{I}_\Gamma(a))=h^0(\mathcal{I}_{\Gamma'}(a))$ for all 
$\Gamma' \subset \Gamma$ such that $|\Gamma'|=|\Gamma|-h^1(\mathcal{I}_\Gamma(a))$. 
\end{prop}
Combining the Cayley-Bacharach Theorem, Proposition~\ref{Kprop} and
Lemma~\ref{Hlemma1} yields our main result:

\begin{thm}\label{MainTh}
  Let $\Gamma \subset \proj{m}$ be a reduced complete intersection of
  hypersurfaces of degrees $d_1, d_2, \ldots, d_m$, and
  let $s = (\sum_{i=1}^m d_i) - m - 1$ as in Theorem~\ref{modernCB}.
  If $1 \leq a \le s$, then the evaluation code $C(\Gamma)_a$ has
  minimum distance $d \ge (\sum_{i=1}^m d_i) - a - (m-1) = s - a + 2$.
\end{thm}

\begin{proof}
  Put $\deg(\Gamma)-i+1 = s-a+2$, so that $i = \deg(\Gamma)-(s-a+1)$.
  Applying Proposition~\ref{Kprop}, we see that the theorem is true
  iff $h^0(\mathcal{I}_{\Gamma'}(a))-h^0(\mathcal{I}_{\Gamma}(a))=0$
  for all subsets $\Gamma'$ with $\deg(\Gamma')=\deg(\Gamma)-(s-a+1)$.
  The modern Cayley-Bacharach Theorem tells us that
  $$h^0(\mathcal{I}_{\Gamma'}(a))-
  h^0(\mathcal{I}_{\Gamma}(a))=h^1(\mathcal{I}_{\Gamma''}(s-a)).$$
  But
  for any subset $\Gamma'' \subset \Gamma$ of $s+1-a$ points,
  Lemma~\ref{Hlemma1} implies that
  $h^1(\mathcal{I}_{\Gamma''}(s-a))=0$.
\end{proof}

\begin{cor}
An evaluation code $C(\Gamma)_a$ obtained from a reduced 
complete intersection $\Gamma$ is MDS iff 
$$h^1(\mathcal{I}_{\Gamma''}(s-a))=0 \mbox{ for all }\Gamma'' \mbox{ such
that }|\Gamma''|=h^1(\mathcal{I}_{\Gamma}(a)).$$
\end{cor}

\begin{proof}
By Proposition~\ref{Kprop}, $C(\Gamma)_a$ is an MDS code iff 
$h^0(\mathcal{I}_\Gamma(a))=h^0(\mathcal{I}_{\Gamma'}(a))$ for all 
$\Gamma' \subset \Gamma$ such that $|\Gamma'|=|\Gamma|-h^1(\mathcal{I}_\Gamma(a))$.
By Cayley-Bacharach,
$h^0(\mathcal{I}_{\Gamma}(a))=h^0(\mathcal{I}_{\Gamma'}(a))$ for all 
subsets $\Gamma'$ of cardinality $i$ iff 
$h^1(\mathcal{I}_{\Gamma-\Gamma'}(s-a))=0$ for 
all subsets $\Gamma'$ of cardinality $i$. Hence, $C(\Gamma)_a$ is 
MDS iff $h^1(\mathcal{I}_{\Gamma''}(s-a))=0$ for all subsets $\Gamma''$
with $|\Gamma''| = |\Gamma|-(|\Gamma| -h^1(\mathcal{I}_\Gamma(a)))$.
\end{proof}
                                                                                
Write $\sigma_\Gamma$ for the largest $i$ such that
$h^1(\mathcal{I}_\Gamma(i)) \ne 0$. A zero-dimensional scheme $\Gamma'$ such 
that
$h^0(\mathcal{I}_\Gamma(\sigma_{\Gamma}))=h^0(\mathcal{I}_{\Gamma'}(\sigma_\Gamma))$ 
for all
$\Gamma' \subset \Gamma$, $|\Gamma'|=|\Gamma|-1$ is called a 
{\it Cayley-Bacharach} scheme.
In \cite{jh}, Hansen showed that if $\Gamma$ is a Cayley-Bacharach scheme,
then $C(\Gamma)_{\sigma_\Gamma}$ is an MDS code. 
Of course, a complete intersection is a Cayley-Bacharach scheme, with
$s = \sigma_\Gamma$, so the complete intersection codes $C(\Gamma)_s$
are MDS. Are there other complete intersection codes which are MDS?
We know that $h^1(\mathcal{I}_{\Gamma''}(s-a))=0$ if $s-a \ge
|\Gamma''|-1$; so we see that a sufficient condition for the MDS property
is $$s-a \ge h^1(\mathcal{I}_{\Gamma}(a))-1.$$ 

\begin{lem}
If $\Gamma$ is a complete intersection, then 
\[
h^1(\mathcal{I}_{\Gamma}(a))=|\Gamma|-h^1(\mathcal{I}_{\Gamma}(s-a)).
\] 
\end{lem}
\begin{proof}
From the four term exact sequence of Section 1.2, it follows that
$h^1(\mathcal{I}_\Gamma(a)) =
|\Gamma|-\dim_\mathbb{K}(R/I_\Gamma)_{a}$.
Thus, it suffices to show
\[
\dim_\mathbb{K}(R/I_\Gamma)_{a} + \dim_\mathbb{K}(R/I_\Gamma)_{s-a} = |\Gamma|
\]
Let $L \in R_1$ be a non-zero divisor on $R/I_\Gamma$ (such an $L$
exists since $R/I_\Gamma$ is Cohen-Macaulay). We pass to the Artinian 
reduction $R/(I_\Gamma+\langle L\rangle)$.  It is easy to see that 
\[
\sum_{i=0}^{s+1} \dim_\mathbb{K}(R/(I_\Gamma+\langle L\rangle))_i = |\Gamma|. 
\]
Since $L$ is not a zero divisor, there is an exact sequence
\[
0 \longrightarrow (R/I_\Gamma)(-1) \stackrel{\cdot L}{\longrightarrow} R/I_\Gamma
 \longrightarrow R/(I_\Gamma+\langle L\rangle) \longrightarrow 0.
\]
From the exact sequence, it follows that
\[
\dim_\mathbb{K}(R/I_\Gamma)_a = \sum_{i=0}^{a}\dim_\mathbb{K}(R/(I_\Gamma+\langle L\rangle)_i).
\]
Similarly, we have
\[
\dim_\mathbb{K}(R/I_\Gamma)_{s-a} = \sum_{i=0}^{s-a}\dim_\mathbb{K}(R/(I_\Gamma+\langle L\rangle)_i).
\]
Now, since $\Gamma$ is a complete intersection, the Hilbert function
of the Artinian reduction is symmetric. So 
\[
\sum_{i=0}^{s-a}\dim_\mathbb{K}(R/(I_\Gamma+\langle L\rangle)_i) = \sum_{i=a+1}^{s+1}\dim_\mathbb{K}(R/(I_\Gamma+\langle L\rangle)_i),
\]
yielding the result.
\end{proof}

Thus, a sufficient condition for the MDS property is that $s-a+1 \ge
\deg(\Gamma)-h^1(\mathcal{I}_{\Gamma}(s-a)) =
\dim_\mathbb{K}(R/I_\Gamma)_{s-a}$.
If $\Gamma$ is a set of collinear points, then 
$\dim_\mathbb{K}(R/I_\Gamma)_{m} = \min\{m+1,|\Gamma| \}$, so a set
of collinear points always gives an MDS code. 

\section{Examples}\label{sec:four}
We now give several examples to illustrate our results. First, 
we quickly review the notation from the previous sections.
We consider codes $C(\Gamma)_a$ constructed by evaluating the homogeneous
polynomials of degree $a$ at the points of a complete intersection 
$\Gamma = X_1\cap \cdots X_m$, where $X_i$ has degree $d_i$.
As in Theorem~\ref{MainTh}, we write $s = (\sum_{i=1}^m d_i) - m - 1$.
Then the result of that theorem says that if $1 \le a \le s$, then 
the minimum distance $d$ of 
the evaluation code satisfies $d \ge s - a + 2$.  

\begin{exm}
Let $x_j$, $0\le j \le m$ be
the homogeneous coordinates on $\proj{m}$, and let
$X_1,\ldots,X_{m-1}$ be the hyperplanes $X_j = V(x_j)$ for $1 \le j
\le m - 1$.  Let $X_m$ be the hypersurface $V(x_m^q - x_0^{q-1}x_m)$.
Then the intersection of the $X_i$ is a complete intersection
$\Gamma$, consisting of the set of affine $\mathbb{F}_q$-rational
points (i.e. points with $x_0 \ne 0$) on the line $L = X_1\cap \cdots
\cap X_{m-1}$.  The evaluation codes in this case are just the usual
extended Reed-Solomon codes, and Theorem \ref{MainTh} yields the
following.  We have $s = m - 1 + q - m - 1 = q - 2$.  If $a \le s$,
then we get that the minimum distance satisfies
$$d \ge q - 2 - a + 2 = q - a = n - k + 1,$$
since the block length
$n$ is $q$, and the dimension $k$ is $a+1$.  Thus we have recovered
the well-known fact that the extended Reed-Solomon codes are MDS
codes.
\end{exm}
\begin{exm}
Second, consider the usual Reed-Muller evaluation codes as in Example
4.5 of \cite{jh1}, where the case $m = 2$ is studied.  The set of all
affine $\mathbb{F}_q$-rational points in $\affine{m}$ is the
projective complete intersection
$$\Gamma = V(x_j^q - x_0^{q-1}x_j : j = 1, \ldots, m).$$
Hence we have $s = mq - m - 1 = m(q - 1) - 1$.
Our Theorem~\ref{MainTh} implies that for the $C(\Gamma)_a$
code with $a \le s$, the minimum distance is bounded below by
$$d \ge s - a + 2 = m(q - 1) - a + 1.$$
We note that this example
shows the type of bound we are considering here is likely to be of
interest in general only when $a$ is relatively large compared to $s$.
For instance, it is known that if $a = \alpha(q - 1) + \beta$, where
$0 \le \beta \le q - 2$, then the exact minimum distance of the
Reed-Muller code is $d = (q - \beta)q^{m - 1 - \alpha}$ (see
\cite{ak}, Corollary 5.5.4, for instance).  If $\alpha < m - 1$, then
our lower bound will be considerably smaller than the actual minimum
distance. On the other hand, if, for example, $a = (m-1)(q-1)$, so
$\alpha = m - 1$ and $\beta = 0$, the actual minimum distance is $d =
q$, while our bound also gives $d \ge m(q-1) + 1 - (m-1)(q-1) = q$.
\end{exm}
\begin{exm}
For our final example, we consider codes related to Hermitian codes.
The evaluation geometric Goppa codes over $\mathbb{F}_{q^2}$ are
defined using the Hermitian curves $X_q = V(x_1^{q+1} - x_2^qx_0 -
x_2x_0^q) \subset \proj{2}$ and the divisors $G = uQ$ where $Q =
[0,0,1]$ is the unique point at infinity on $X_q$. There are precisely
$q^3$ affine $\mathbb{F}_{q^2}$-rational points on $X_q$.  However the
$\Gamma$ consisting of all of them is not a projective complete
intersection. To construct codes for which our main results apply, we
let
$$F(x_0,x_1,x_2) = \prod_{\{\alpha \in \mathbb{F}_{q^2} :
  \alpha^q+\alpha\ne 0\}}(x_2-\alpha x_0).$$
Then $\Gamma = X_q \cap
V(F)$ consists of the $q^3 - q$ $\mathbb{F}_{q^2}$-rational points on
$X_q$ with $x_1 \ne 0$ (all in the affine part of the plane).  In a
very precise sense (see \cite{lsh}), the evaluation codes
$C(\Gamma)_a$ are related to the usual Hermitian codes constructed
using the divisor $D$ consisting of all $\mathbb{F}_{q^2}$-rational
points in the same way that Reed-Solomon codes are related to the
extended Reed-Solomon codes.

As in the Reed-Muller case, our bound only gives sharp results when
the degree $a$ is large relative to $s$.  Since the equations defining
$\Gamma$ have degrees $d_1 = q+1$ and $d_2 = q^2 - q$, we have $s =
q^2 - 2$.  For example, with $a = q^2 - q$, our Theorem~\ref{MainTh}
yields $d \ge s + 2 - a = q$.  By way of comparison, the usual
Hermitian evaluation code constructed using $L(uQ)$ for $u = a(q+1) =
q^3 - q$ (the maximum pole order at $Q$ of the functions corresponding
to the elements of $R_a$) also has $d = q^3 - (q^3 - q) = q$ by
\cite{sti}, Proposition VII.4.3.  Note that our code has block length
$n = q^3 - q$ rather than $q^3$, and the dimension is also one
less than the dimension of the corresponding usual Hermitian code
because the polynomial $F$ has degree $a = q^2 - q$.
\end{exm}

There is an extension of the notion of a residual scheme from the case
when $\Gamma$ is a complete intersection to the case when $\Gamma$ is
arithmetically Gorenstein. It seems reasonable to expect that similar
methods would yield bounds on the minimum distance in this case; we
hope to study this question in a future paper. We note that in
\cite{ep}, Eisenbud and Popescu use the (local) Gorenstein property to
give a proof of Goppa Duality.

\bigskip
\noindent{\bf Acknowledgment} 
This collaboration began while the authors were members of MSRI; we
also thank the Institute for Scientific Computation at Texas A\&M for
providing logistical support.

\renewcommand{\baselinestretch}{1.0}
\small\normalsize 

\bibliographystyle{amsalpha}

\end{document}